\documentstyle[11pt]{article}

\newenvironment{ab}{\begin{quote}\begin{scriptsize}\begin{it}
{\bf Abstract}\,--\,}{\end{it}\end{scriptsize}\end{quote}}

\newenvironment{re}{\begin{quote}\begin{scriptsize}\begin{it}
{\bf R\'esum\'e}\,--\,}{\end{it}\end{scriptsize}\end{quote}}

\newtheorem{thm}{Theorem}[section]

\newcommand{\C}{\mbox{$\bf C$}}     
\newcommand{\Z}{\mbox{$\bf Z$}}     

\newcommand{\ssp}[1]{\mbox{$\scriptscriptstyle {#1}$}}

\newcommand{\im}{{\rm im}\,}        
\newcommand{\coker}{{\rm coker}\,}  

\newcommand{\Pic}{{\rm Pic}\,}      
\newcommand{\bPic}{\mbox{$\bf Pic$}} 
\newcommand{\Spec}{{\rm Spec}\,}    

\newcommand{\G}{\mbox{$\bf G$}}     
\newcommand{\Alb}{{\rm Alb}\,}

\newcommand{\by}[1]{\stackrel{#1}{\rightarrow}}

\renewcommand{\tilde}{\widetilde}

\newcommand{\df}{\mbox{\,$\stackrel{\ssp{\rm def}}{=}$}\,}
\newcommand{\ie}{{\it i.e.\/},\ }
\newcommand{\cf}{{\it c.f.\/}\ }
\newcommand{\eg}{{\it e.g.\/},\  }

\renewcommand{\bar}{\overline}

\newcommand{\longto}{\longrightarrow}

\newcommand{\boxtensor}{{\Box\kern-9.03pt\raise1.42pt\hbox{$\times$}}}
\newcommand{\supp}{{\rm supp}\,}
\renewcommand{\d}{\mbox{\LARGE $\cdot $}}

\newcommand{\Div}{{\rm Div}\,}
\newcommand{\arc}[1]{{#1}\kern-16pt\raise7pt\hbox{\Large{$\frown$}}}

\newcommand{\longinto}{\lhook\joinrel\kern-3pt\hbox to
100pt{\rightarrowfill}}           

\newcommand{\propsubset}
{\mbox{$\textstyle{
\subseteq_{\kern-5pt\raise-1pt\hbox{\mbox{\tiny{$/$}}}}}$}}

\newcounter{elno}                

\newcounter{example}[section]
\def\theexample{\thesection.\arabic{example}}

\newcounter{exercise}[section]
\def\theexercise{\thesection.\arabic{exercise}}

\newcommand{\ccL}{{\cal L}}

\newcommand{\cO}{{\cal O}}

\begin{document}
G\'eom\'etrie alg\'ebrique/{\it Algebraic Geometry}\\[4pt]
\begin{center}
{\large\bf Albanese and Picard 1-motives}\\[4pt]
Luca {\sc Barbieri Viale} and Vasudevan {\sc
Srinivas}\\[4pt]
\end{center}
\begin{ab}
We define, in a purely algebraic way, 1-motives
$\Alb^{+}(X)$,  $\Alb^{-}(X)$, $\Pic^{+}(X)$ and
$\Pic^{-}(X)$ associated with any algebraic scheme $X$
over an algebraically closed field of characteristic zero.
For $X$ over $\C$ of dimension $n$ the Hodge realizations are,
respectively, $H^{2n-1}(X,\Z (n))/({\rm torsion})$,
$H_{1}(X,\Z)/({\rm torsion})$, $H^{1}(X,\Z (1))$ and
$H_{2n-1}(X,\Z (1-n))/({\rm torsion})$.
\end{ab}
\begin{center}
{\bf 1-motifs de Albanese et de Picard}
\end{center}
\begin{re}
Nous d\'efinissons, de fa\c con purement alg\'ebrique,
les 1-motifs $\Alb^{+}(X)$,  $\Alb^{-}(X)$,
$\Pic^{+}(X)$ et $\Pic^{-}(X)$ associ\'es au tout
sch\'ema alg\'ebrique $X$ sur un corps de
charact\'eristique 0. Pour $X$ sur $\C$ les r\'ealisations
de Hodge sont, respectivement, $H^{2n-1}(X,\Z (n))/({\rm
torsion})$, $H_{1}(X,\Z)/({\rm torsion})$, $H^{1}(X,\Z
(1))$ et $H_{2n-1}(X,\Z (1-n))/({\rm torsion})$.
\end{re}
\vspace{0.5cm}
\begin{small}
{\bf Version fran\c caise abr\'eg\'ee}\,--\, Dans cette
Note nous annon\c cons certains r\'esultats dont les preuves
se trouvent dans [1]. Nous donnons une construction
alg\'ebrique pour les 1-motifs ``d'Albanese et de Picard'' en
accord avec la conjecture de Deligne selon laquelle les
1-motifs associ\'es \`a la structure de Hodge mixte sur la
cohomologie des variet\'es alg\'ebriques admettent une
definition alg\'ebrique (voir  [3], \S 10.4.1). Etant
donn\'e que les structures de Hodge mixtes  $H^{2n-1}(X,\Z
(n))/{\rm (torsion)}$, $H_{1}(X,\Z)/({\rm torsion})$,
$H^{1}(X,\Z (1))$ et $H_{2n-1}(X,\Z (1-n))/({\rm torsion})$
pour $X$ sur $\C$ de dimension $n$, sont de niveau $\leq 1$,
elles sont conjecturellement alg\'ebrisables. Soit $X$ un sch\'ema
alg\'ebrique sur un corps $k$ alg\'ebriquement clos de
charact\'eristique 0. Soit $f:\tilde X\to X$ une r\'esolution des
singularit\'es,  $\bar X$ une compactification lisse de
$\tilde X$, $Y$ le diviseur \`a croisements normaux et $Y=\cup Y_i$
la r\'eunion des composantes lisses. Soit $\Pic (\bar X, Y)$ le
groupe des classes d'isomorphisme des pairs $(\ccL, \varphi)$,
telles que $\ccL$ soit un fibr\'e en droites sur $\bar X$ et
$\varphi :\ccL\mid_Y\cong \cO_Y$ une trivialisation sur $Y$.\\[2pt]
$\ominus$ {\sc Lemme}\,--\, {\it Le foncteur de Picard rigidifi\'e
$T\leadsto  \Pic (\bar X\times_k T, Y\times_k T)$ est repr\'esentable.
La composante connexe de l'identit\'e  $\Pic^0_{(\bar X, Y)/k}$ est un
sch\'ema semi-ab\'elien, repr\'esentable par une extension
$$0\to T(\bar{X},Y)\to
\Pic^0_{(\bar X, Y)/k}\to [\ker^0 \Pic^0_{{\bar X}/k}\to
\oplus\Pic^0_{Y_i/k}]\to 0$$  d'un sch\'ema ab\'elien par un
tore.}\\[2pt] Pour chaque $D$ disjoint de $Y$ il y a une pair
canonique $(\cO(D),1)$ dans $\Pic (\bar X, Y)$. Soit $S$ le lieu
singulier de $X$ et soit $\bar S$ la fermeture de $f^{-1}(S)$ dans
$\bar X$. Alors $\Div_{\bar S}^0(\bar X,Y)$ est le groupe des
diviseurs engendr\'e par les composantes propres de $f^{-1}(S)$ qui
sont alg\'ebriquement nulles. Soit $\Div_{\bar S/S}^0(\bar X,Y)$ le
sous-groupe des diviseurs ayant image directe nulle sur $S$. Nous
d\'efinissons le 1-motif ``Picard homologique'' de la fa\c con
suivante: $$\Pic^-(X)\df [\Div_{\bar S/S}^0(\bar X,Y)\to\Pic^0(\bar
X,Y)]$$ Le 1-motif ``Albanese cohomologique'' $\Alb^+(X)$ est, par
d\'efinition, le dual de Cartier de $\Pic^-(X)$.\\[4pt]  $\ominus$
{\sc Th\'eor\`eme}\,--\,  {\it Si $X$ est d\'efini sur $\C$ et
$n=\dim X$, alors la r\'ealisation de Hodge de $\Pic^-(X)$ est
$H_{2n-1}(X,\Z(1-n))/{\rm (torsion)}$.}\\[4pt] La
r\'ealisation de Hodge de $\Alb^+(X)$ est alors
$H^{2n-1}(X,\Z(n))/{\rm (torsion)}$; par example,
$\Alb^+(X)\cong J^n(X)$ lorsque $X$ est une vari\'et\'e
projective complexe. Soit $f_{\d}: X_{\d}\to X$ un
hyperrecouvrement de $X$; alors $X_{\d}$ est un sch\'ema
simplicial avec composantes lisses. Nous considerons une
compactification propre et lisse $\bar X_{\d}$ \`a
croisements normaux $Y_{\d}$; nous d\'esignons par
$\Div_{Y_{\d}}(\bar X_{\d})$ le noyau de $d_0^*-d_1^*: \Div_{Y_0}(\bar
X_0)\to\Div_{Y_1}(\bar X_1)$\\[2pt]   $\oplus$ {\sc
Lemme}\,--\, {\it Le $fpqc$-faisceau $\bPic_{{\bar
X_{\d}}/k}$ associ\'e au foncteur de Picard simplicial
$T\leadsto \bPic (\bar X_{\d}\times_k T)$ est
repr\'esentable. La composante connexe de l'identit\'e est un
sch\'ema semi-ab\'elien, repr\'esentable par une extension $$0\to
T(\bar X_{\d})\to\bPic^0_{{\bar X_{\d}}/k} \to [\ker^0\Pic_{{\bar
X_0}/k}^0\to\Pic_{{\bar X_1}/k}^0] \to 0$$
d'un sch\'ema ab\'elien par un tore.}\\[4pt] Nous d\'esignons par
$\Div_{Y_{\d}}^0(\bar X_{\d})$ le sous-groupe des diviseurs ayant
image dans $\bPic^0$. Nous d\'efinissons le 1-motif ``Picard
cohomologique'' de $X$ $$\Pic^+(X)\df [\Div_{Y_{\d}}^0(\bar
X_{\d})\to \bPic^0(\bar X_{\d})]$$ Le dual de Cartier de $\Pic^+(X)$
est le 1-motif ``Albanese homologique'' $\Alb^-(X)$.\\[2pt] $\oplus$
{\sc Th\'eor\`eme}\,--\, {\it Si $X$ est defini sur $\C$, alors la
r\'ealisation de Hodge de $\Pic^+(X)$ est $H^{1}(X,\Z(1))$.}\\[4pt]
La r\'ealisation de Hodge de $\Alb^-(X)$ est $H_{1}(X,\Z)/{\rm
(torsion)}$. Certaines applications \`a l'\'etude des cycles
alg\'ebriques ainsi que d'autres propri\'et\'es concernant les
1-motifs d\'efinis ci-dessus para\^\i trons dans
[1].\\[2pt]\end{small} \footnoterule\hfill\\[4pt]
In this Note we announce some
results whose proof will appear in [1]. We give an algebraic
construction for ``Albanese and Picard'' 1-motives, in accordance
with Deligne's conjecture that the 1-motives associated to mixed
Hodge structures on the cohomology groups of an algebraic
variety should be algebraically defined (see [3], \S 10.4.1).
The mixed Hodge structures $H^{2n-1}(X,\Z (n))/{\rm
(torsion)}$, $H_{1}(X,\Z)/({\rm torsion})$, $H^{1}(X,\Z (1))$
and $H_{2n-1}(X,\Z (1-n))/({\rm torsion})$ for $X$ over $\C$
of dimension $n$ are of the appropriate type, \ie are of niveau
$\leq 1$, and hence the associated 1-motives conjecturally do have
an algebraic description.

Let $X$ be an algebraic scheme over an algebraically closed field
$k$ of characteristic zero.\\[4pt]
{\sc Homological Picard and cohomological Albanese 1-motives}\,--\,
Let $f:\tilde X\to X$ be a resolution of singularities and let $\bar
X$ be a smooth compactification of $\tilde X$ with normal crossing
boundary $Y$, which has smooth irreducible components $Y_i$. Let
$\pi_{\bar{X}}:\bar{X}\to\Spec k$, $\pi_Y:Y\to\Spec k$ be the structure
morphisms. Then the $fpqc$-sheaves $(\pi_{\bar{X}})_*\G_{m,\bar{X}}$ and
$(\pi_Y)_*\G_{m,Y}$ are represented by algebraic $k$-tori, as is
$T(\bar{X},Y)=\coker\left((\pi_{\bar{X}})_*\G_{m,\bar{X}} \to
(\pi_Y)_*\G_{m,Y}\right)$. If $G\to H$ is a homomorhism of $k$-group
schemes which are locally of finite type over $k$, let $\ker^0 G\to H$
denote the connected component of the identity of the $k$-group scheme
$\ker G\to H$. \\[4pt]
$\ominus$ {\sc Lemma}\,--\, {\it The rigidified Picard functor
$T\leadsto  \Pic (\bar X\times_k T, Y\times_k
T)$ is representable by a $k$-group scheme which is locally of finite
type, and whose group of $k$-points is $\Pic(\bar X,Y)$. The connected
component of the identity $\Pic^0_{(\bar X,
Y)/k}$ is a semi-abelian scheme which can be represented as an extension
$$0\to T(\bar{X},Y)\to
\Pic^0_{(\bar X, Y)/k}\to [\ker^0 \Pic^0_{{\bar X}/k}\to
\oplus\Pic^0_{Y_i/k}]\to 0$$
of an abelian scheme by a torus.}\\[4pt]
We let $\Pic^0(\bar X,Y)\subset \Pic(\bar X,Y)$ denote the subgroup of
$k$-points of $\Pic^0_{(\bar X,Y)/k}$. Any divisor $D$ on $\bar X$ with
support disjoint from $Y$ has a class $(\cO_{\bar X}(D),1)\in\Pic(\bar
X,Y)$; we say $D$ is algebraically equivalent to 0 relative to $Y$ if
this class is in the subgroup $\Pic^0(\bar X,Y)$.
Let $S$ denote the singular locus of $X$, and $\bar{S}$ the closure in
$\bar{X}$ of $f^{-1}(S)$. Let $\Div^0_{\bar{S}}(\bar X,Y)$ denote the
group of divisors $D$ on $\bar X$ such that $\supp D$ is disjoint from
$Y$, and contained in $\bar{S}$, and also $D$ is algebraically equivalent
to 0 relative to $Y$. Let $\Div^0_{\bar S/S}(\bar X,Y)\subset \Div^0_{\bar
S}(\bar X,Y)$ denote the subgroup of such divisors $D$ which have
vanishing push-forward under $f$.
We can now define the ``homological Picard'' 1-motive as
follows:
$$\Pic^-(X)\df [\Div_{\bar S/S}^0(\bar X,Y)\to\Pic^0(\bar X,Y).]$$
The ``cohomological Albanese'' 1-motive $\Alb^+(X)$ is
defined to be the Cartier dual of $\Pic^-(X)$.\\[4pt]
$\ominus$ {\sc Theorem}\,--\,  {\it If $X$ is defined over
$\C$ and $n=\dim X$, then the Hodge realization of
$\Pic^-(X)$ is $H_{2n-1}(X,\Z(1-n))/{\rm
(torsion)}$.}\\[4pt] By duality, we have that the Hodge realization of
$\Alb^+(X)$ is $H^{2n-1}(X,\Z(n))/{\rm (torsion)}$, \eg
$\Alb^+(X)\cong J^n(X)$ if $X$ is a projective complex
variety, where $J^n(X)$ is the semi-abelain variety considered in
[2].\\[4pt]
{\sc Cohomological Picard and homological Albanese 1-motives}\,--\,
We let $f_{\d}: X_{\d}\to X$ be a proper smooth
hypercovering of $X$; then $X_{\d}$ is a simplicial $k$-scheme
with smooth components. Consider a proper smooth
compactification  $\bar X_{\d}$ with normal crossing boundary
$Y_{\d}$, and denote by $\Div_{Y_{\d}}(\bar X_{\d})$
the group of divisors on $X_0$ which are supported on $Y_0$, and have
zero pullback on $X_1$, \ie by definition $$\Div_{Y_{\d}}(\bar
X_{\d})\df \ker\Div_{Y_0}(\bar
X_0)\by{d_0^*-d_1^*}\Div_{Y_1}(\bar X_1).$$
Let $\pi_i:\bar X_i\to\Spec k$ be the structure morphism. Let
$T(\bar  X_{\d})$ be the $k$-torus representing the $fpqc$-sheaf
\[\frac{\ker [(\pi_1)_*\G_{m,\bar
X_1}\to (\pi_2)_*\G_{m,\bar X_2}]}{\im [(\pi_0)_*\G_{m,\bar X_0}\to
(\pi_1)_*\G_{m,\bar X_1}]}.\]
$\oplus$ {\sc Lemma}\,--\, {\it The sheaf $\bPic_{{\bar
X_{\d}}/k}$, with respect to $fpqc$-topology, associated
to the simplicial Picard functor $T\leadsto \bPic
(\bar X_{\d}\times_k T)$, is representable by a group scheme locally of
finite type over $k$, with $k$-points $\Pic(\bar X_{\d})$. The connected
component of the identity is a semi-abelian scheme, which can be
represented as an extension
$$0\to T(\bar X_{\d})\to\bPic^0_{{\bar X_{\d}}/k} \to
[\ker^0 \Pic_{{\bar X_0}/k}^0\to\Pic_{{\bar X_1}/k}^0] \to 0$$
of an abelian scheme by a torus.}\\[4pt]
Let $\Div_{Y_{\d}}^0(\bar X_{\d})$ denote the
subgroup of $\Div_{Y_{\d}}(\bar X_{\d})$ of those divisors which are
mapped to $\bPic^0(k)$ under the canonical mapping. We then define the
``cohomological Picard'' 1-motive of $X$ as follows:
$$\Pic^+(X)\df [\Div_{Y_{\d}}^0(\bar X_{\d})\to \bPic^0(\bar X_{\d})].$$
The ``homological Albanese'' 1-motive $\Alb^-(X)$ is defined to be the
Cartier dual of $\Pic^+(X)$.\\[4pt]
$\oplus$ {\sc Theorem}\,--\, {\it If $X$ is defined over $\C$, then the
Hodge realization of $\Pic^+(X)$ is $H^{1}(X,\Z(1))$.}\\[4pt]
We then also have that the Hodge realization of $\Alb^-(X)$ is
$H_{1}(X,\Z)/{\rm (torsion)}$.\\[4pt]
{\sc Remarks}\,--\, We keep the same notation as above.
\begin{enumerate}
\item We can show that our definitions are independent of
choices of resolutions, compactifications or hypercoverings.
\item If $X$ is normal and projective over $k= \bar k$ then $\Alb^+(X)$
is an abelian variety, which may be identified with $\Alb (\tilde X)$,
the Albanese variety of a resolution. If $X$ is projective, but not
necessarily normal, then $\Alb^+(X)$ is a semi-abelian variety. If
$X$ is smooth, possibly open, then $\Alb^-(X)$ is a semi-abelian
variety.
\item Let $X$ be irreducible and projective of dimension $n$, and let
$X_{\rm reg}$ denote the set of smooth points of $X$; we may also
regard $X_{\rm reg}$ as an open subscheme of any given resolution of
singularities $\tilde{X}$.  Let $x_0\in X_{\rm reg}$ be a smooth
point and let $a_{x_0}:\tilde{X}\to \Alb(\tilde{X})$ be the
corresponding Albanese mapping;  considering $\Alb^+$ as a torus
bundle over $\Alb(\tilde{X})$, and pulling back along $a_{x_0}$, one
can see it has a natural trivialization on $X_{\rm reg}$, yielding a
section $a_{x_0}^+: X_{\rm reg}\to \Alb^+(X)$. This map yields a
universal regular homomorphism $CH^{n}(X)_{\deg 0}\to \Alb^+(X)$
from the ``cohomological'' Levine-Weibel Chow group of zero cycles
to semi-abelian varieties, and the corresponding Roitman torsion
theorem can be proved (\cf [2]). \item We clearly have a
commutatitve square of 1-motives $$\begin{array}{ccc} \Pic^-(\tilde
X)&\longto &\Pic^+(\tilde X) \\ \downarrow & &\downarrow \\
\Pic^-(X) & \longto & \Pic^+(X_{\rm reg})\\ \end{array}$$
In particular, by taking Cartier duals we have a map
$\Alb^-(X_{\rm reg})\to \Alb^+(X)$. Then there is a morphism
$X_{\rm reg}\to  \Alb^-(X_{\rm reg})$ which lifts the section
$a_{x_0}^+$  above when $X$ is proper. In fact, $\Alb^-$ of a smooth
variety is the  universal semi-abelian variety associated to it by
Serre [6], \S 5.  \item The claimed representabilities of
``rigidified'' and ``simplicial''  Picard functors can be obtained
{\it via}\, the classical theorems due to Grothendieck [4] and Murre
[5].   \item For further properties we refer to [1] where we will
consider  $\ell$-adic and De Rham realizations; functoriality
properties, algebraic Gysin maps and Lefschetz theorems (\cf [2])
are also expected. There are also analogous definitions of the
1-motives in the case when the ground field $k$ is an arbitrary
field of characteristic 0. However, the tori appearing need not be
split, for example. The general case is reduced to the case of an
algebraically closed field using Galois descent.  \end{enumerate}

\begin{small}
Ce travail a \'et\'e r\'ealis\'e lors de s\'ejours des
auteurs \`a Paris VII, TIFR Mumbai, DIMA Genova et
ICTP Trieste, que nous remercions.\\[4pt]
{\sc R\'ef\'erences}

[1] {\sc L. Barbieri Viale} et {\sc V. Srinivas}:
Albanese and Picard 1-motives, travail en pr\'eparation.

[2] {\sc J. Biswas} et {\sc V. Srinivas}: Roitman's
theorem  for singular projective varieties, preprint, 1996.

[3] {\sc P. Deligne}: Th\'eorie de Hodge III, {\it Publ.
Math.}\, IHES {\bf 44} (1974) 5--78.

[4] {\sc A. Grothendieck}: Technique de descente et
th\'eor\`emes d'existence en g\'eom\'etrie alg\'ebrique
V-VI Les sch\'emas de Picard, dans ``Fondements de la
G\'eom\'etrie Alg\'ebrique'' Extraits du {\it  S\'eminaire}\,
{\sc Bourbaki} 1957/62.

[5] {\sc J. Murre}: On contravariant functors from the
category of preschemes over a field into the category of
abelian groups (with application to the Picard functor),
{\it Publ. Math.}\, IHES {\bf 23} (1964) 581--619.

[6] {\sc J.P. Serre}: Morphismes universels et
vari\'et\'es d'Albanese, dans ``Vari\'et\'es de Picard'' ENS
{\it S\'eminaire}\, {\sc C. Chevalley} 3e ann\'ee: 1958/59.\\

\begin{flushright}
\footnoterule

{\it L. B. V.: Dipartimento di Matematica, Universit\`a di
Genova\\ Via Dodecaneso, 35, I-16146 Genova,
Italia;

V.S.: Tata Institute of Fundamental Research, School
of Mathematics\\ Homi Bhabha Road, 400 005 Mumbai, India.}
\end{flushright}
\end{small}
\end{document}